%% file: jugprob.tex
\documentclass[letterpaper]{amsart}
\usepackage{amsmath,amsthm,amssymb}
\usepackage{tabularx,epsfig,graphicx}
\usepackage{eepic}
\usepackage{psfrag}
\usepackage{url}
\usepackage{mdwlist}
\usepackage{fullpage}

\begin{document}

\renewcommand{\labelenumi}{\theenumi.}

\bibliographystyle{plain}

\input{macros.tex}
\input{jugmac.tex}

\title{Juggling Probabilities}

\date{\today}


\author{Gregory S. Warrington}



%
%

\maketitle



\section{Introduction}
\label{sec:intro}

Imagine\footnote{The author would like to thank Peter Doyle and J.
  Laurie Snell for their substantial help in simplifying the original
  proof of the main result of this paper.}  yourself effortlessly
juggling five balls in a high, lazy pattern.  Your right hand catches a
ball and immediately throws it.  One second later, it is your left
hand's turn to catch and throw a ball.  Then it is your right hand's
turn again$\fourdots$ Some of your throws may be low; some high.  Some
balls go straight up; others cross over to the opposite hand.  At most
times, one ball lands; occasionally a hand remains empty.  But the
alternating cadence of your hands is unwavering.

Suppose that, while you are juggling, we momentarily pause time.
Certain balls are in the air --- you have already thrown them.  To
avoid dropping the ball, you must make catches at certain times in the
future.  If your previous few throws were all low, perhaps you are
only committed at 1,2,4,6 and 7 seconds in the future.  On the other
hand, if you had just vigorously launched a ball, you might be
committed at 1,2,3,4 and 10 seconds in the future.  The set of
``committed times'' is your \textit{landing state} (``state'' for
short).  As you juggle, you wander from state to state according to
what throw you have most recently made.

Our goal in this paper is to answer the following
\newtheorem*{quest}{Random Question}
\begin{quest}\label{mainquest}
  What fraction of the time is spent in any given state if every throw
  is chosen randomly?
\end{quest}
The answer, of course, depends on how we specify the following parameters:
\begin{enumerate}
\item the possible states,
\item the legal throws from a given state and
\item the probability of making each legal throw.
\end{enumerate}
There are countless ways to make the above specifications, but we will
begin our investigations with the model that most closely mirrors what
people envision as juggling.  In the next section, we describe the
juggling universe determined by this standard model and answer our
Random Question.  In \secref{sec:main} we generalize the notion of a
state in order to give a simple proof of our answer.  Finally, in
\secref{sec:models}, we describe variations of our model.

\section{Juggling states}
\label{sec:prelim}

\begin{align*}
\text{Jugglers have been known to} &\ \text{\{toss, roll, spin, bounce, drop, kick\}}\\
  \text{their} & \ \text{\{balls, rings, clubs, flames, chairs\}}\\   
  \text{while accompanying their act with} & \ \text{\{patter, music, dance, flourishes\}}.
\end{align*}
We, however, will strip juggling down to what is (arguably) its
essentials.  Our juggler will juggle only balls and only in the air.
He will not attempt jokes.  Whether a hand is under a leg or behind
the back when making a catch will not be noted.  \emph{We only record
  the order in which the balls are thrown and caught}.  While clearly
ignoring much of what makes juggling visually intriguing, our
single-mindedness will free us to highlight the inherently mathematical
nature of juggling.

Our juggler will:
\begin{enumerate}
\item throw alternately from each hand, once a second,
\item throw a ball instantaneously upon catching it, 
\item be assumed to have always been (and be forevermore) juggling\\
  (this allows us to avoid starting/stopping issues),
\item be named \magnus\footnote{After Magnus Nicholls, a juggler
    prominent during the early 1900s.  He is reputed to be the first
    person to juggle five clubs.  Reliable details of his life are
    outclassed by fascinating stories.}.
\end{enumerate}

Every ball \magnuss throws lands some (integral) number of seconds in
the future.  By abuse of physics, we will refer to this parameter as
the \textit{height} of the throw.  Until mentioned otherwise, we now
fix some maximum height $h$ above which \magnuss is too weak to throw.

In the introduction, we referred to a state as a set of future times
at which balls are landing.  It will be more convenient throughout the
paper to define a \emph{(landing) state} $\stv$ as an $h$-tuple in the
alphabet $\{\cat,\ \emp\}$.  Given a state $\stv$, we will write
$\nu_t$ for its $t$-th component and write $\stv$ as 
\begin{equation*}
  \stv = \text{\fbox{\kern-3pt \begin{tabular}{ZZZZ} 
        $\snuo$ & $\snut$ & $\cdots$ & $\snuh$
  \end{tabular}}}\,.
\end{equation*}
If $\nu_t = \cat$, then \magnuss will catch a ball $t$ seconds in the
future.  If $\nu_t = \emp$, then no ball lands $t$ seconds in the
future.  We make that convention that $\nu_t = \emp$ for $t > h$.  To
recover our original notion of a state from $\stv$, take the set $\{t:
\nu_t = \cat\}$.  We will denote the set of states that are $h$-tuples
by $\stn$.  The subset $\stnk\subset \stn$ will contain only those
states with $f$ $\emp$'s.

As an example, consider the state $\stv = \xxoxoxxxoo\in\stp{10}$.
Assuming that \magnuss has just thrown out of his left hand, at 1 and
7 seconds in the future he will be (instantaneously) throwing/catching
out of his right hand; at 2,4,6 and 8 seconds in future he will be
(instantaneously) throwing/catching out of his left hand.

When the first element of our state $\stv$ is $\emp$, we will often
refer to what \magnuss does during the next second as making a
``height 0'' throw.  Upon making a throw (even a 0), \magnuss commits
himself to a (usually different) new state $\stw$:
\begin{enumerate}
\item If the throw is a 0, $\stw = \fbox{\kern-3pt
    \begin{tabular}{ZZZZZ} $\snut$ & $\snur$ & $\cdots$ & $\snuh$ &
      $\emp$ \end{tabular}}$\,.
\item If the throw is of height $t \geq 1$, then
  $\stw = \fbox{\kern-3pt
    \begin{tabular}{ZZZZZZZ} $\snut$ & $\snur$ & $\cdots$ & $\cat$ & $\cdots$ & $\snuh$ &  $\emp$ \end{tabular}}$\,, where the $\cat$ is placed in the $t$-th position (replacing $\nu_{t+1}$).
\end{enumerate}

As an example, below we illustrate successive throws of height 2, 0
and 6 from an initial state in which balls are landing at 1, 4 and 6
seconds in the future:
\begin{equation*}
  \xooxox \starr{2} \oxxoxo \starr{0} \xxoxoo \starr{6} \xoxoox.
\end{equation*}
In \figref{fig:g34}, we show all possible states and transitions when
\magnuss has 3 balls and is not able to throw anything higher than a
4.

\mymslfig{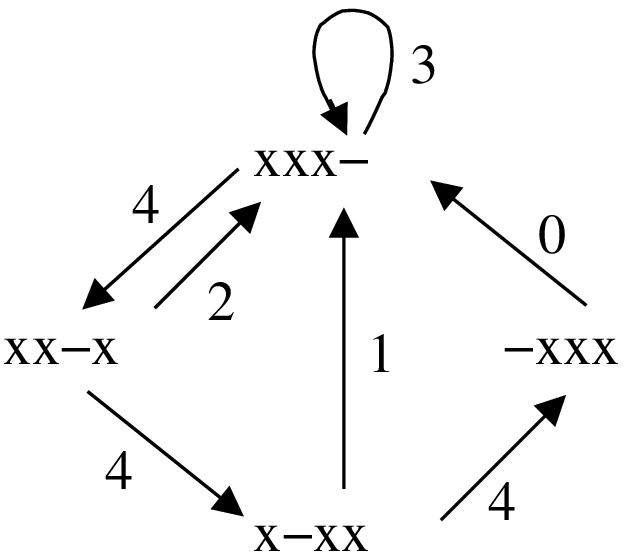}{Juggling 3 balls with maximum height 4.  The edges are
  labeled with the throw height.}

\begin{remark}
  In this section, we have introduced the fundamentals of ``siteswap''
  notation --- albeit in an untraditional form.  As with many great
  ideas, the provenance of siteswap notation is somewhat murky.
  However, we can safely say that this notation was introduced
  (independently, and in various forms) in the mid 1980s by Paul
  Klimek; Colin Wright and members of the Cambridge Juggling Club; and
  Bengt Magnusson and Bruce Tiemann.  For a more leisurely
  introduction and a bibliography of papers relating to siteswap
  notation, the reader should consult \cite{BEGW}.  General internet
  resources on juggling include \cite{juggen} and \cite{jugsoft}.  A
  recent book devoted to the mathematics of juggling is
  \cite{Polster}.
  
  Both as a useful notation for juggling and for interesting
  enumerative combinatorics, it is preferable to introduce siteswap
  notation by defining patterns as repeating sequences of throws.
  Such patterns correspond to cycles in certain graphs; such a graph
  can be found in \figref{fig:g34}.  Many computer programs
  have been written to animate these patterns.  See \cite{magnus} for
  one cowritten by the author, or \cite{jugsoft} for a comprehensive
  list.
\end{remark}


To answer our Random Question, we need an explicit description of the
set of possible states \magnuss can pass through and of the possible
transitions among them.  To describe these transitions, we introduce a
family of \emph{throwing }operators $\adv{j}$ acting on these states.
Intuitively, $\adv{0}$ corresponds to \magnuss waiting when no ball
lands and $\adv{j}$ corresponds to \magnuss making a height $j$ throw
when a ball does land.  Precisely,
\begin{equation}
  \adv{j}(\stv) = 
  \begin{cases}
    \shenuh, &\text{ if } j = 0,\phantom{\sumsb{M\\M}}\\
    \insjq, &\text{ if } j > 0.
  \end{cases}
\end{equation}

For the remainder of this and the next section, we will fix $h > 0$
and $0 \leq f \leq h$.  We are now ready to define the directed graph
which specifies the possible states for \magnuss along with his
possible transitions.  Informally, the valid transitions correspond
either to waiting or to making throws to empty landing times.
\begin{dfn}
  The \emph{state graph} $\grastd{h}{f}$ is the directed graph with
  \begin{enumerate}
  \item Vertices indexed by $\stnk$.
  \item An edge from $\stv$ to $\stw$ whenever 
    \begin{enumerate}
    \item $\stw = \adv{j}(\stv)$ for some $j$ with $1 \leq j \leq h$,
      $\stvr_1 = \cat$ and $\nu_{j+1} = \emp$; or
    \item $\stw = \adv{0}(\stv)$ and $\stvr_1 = \emp$.
    \end{enumerate}
    The set of all edges is denoted $\edges(\grastd{h}{f})$.
  \end{enumerate}
  If $(\stv,\stw)\in\edges(\grastd{h}{f})$, then we will refer to $\stv$ as
  a \emph{precursor} of $\stw$.
\end{dfn}
\figref{fig:g35} illustrates $\grastd{5}{2}$.
\mymslfig{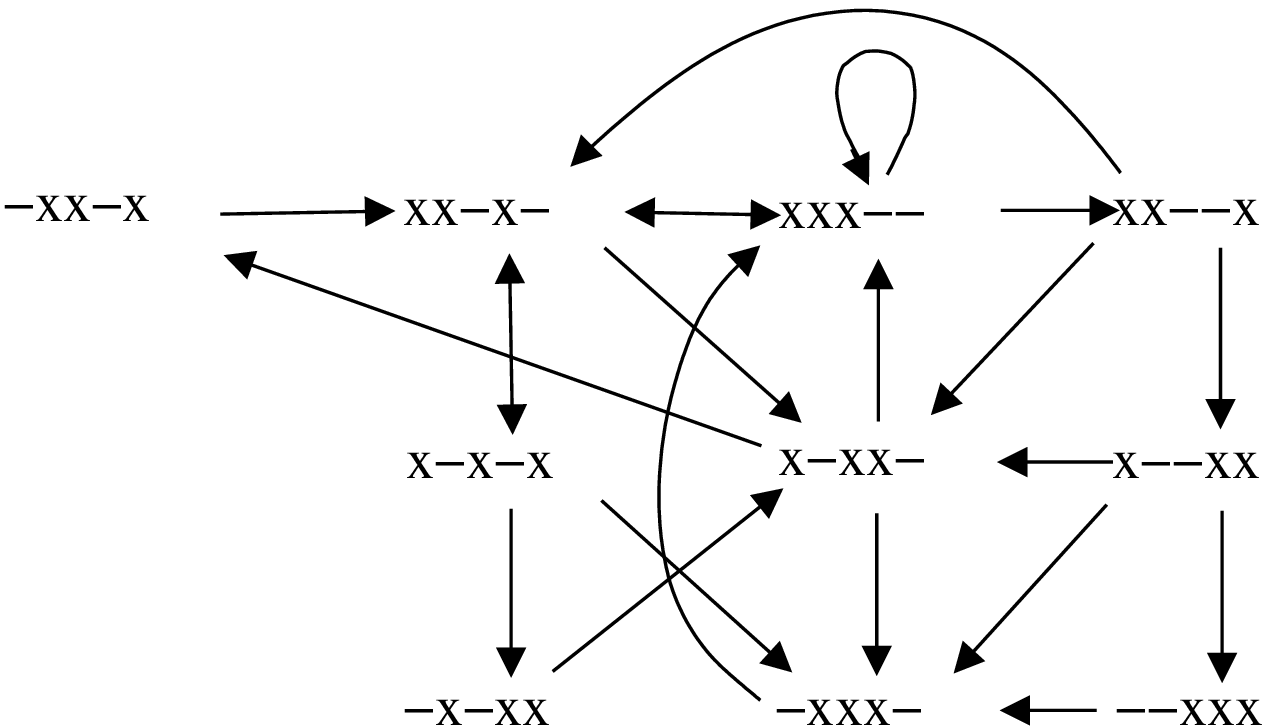}{}

Let us revisit the scenario introduced at the beginning of the paper
(this time with \magnuss juggling).  When we freeze the action,
\magnuss is in a particular state --- he is committed to making
certain catches.  If no ball is landing in one second, he is forced to
wait.  If a ball is landing, however, he can throw to any height
($\leq h$) that will not lead to two balls landing at the same time.
Our assumption in the following is that \magnuss chooses, with equal
probability, one of these throws.  In this manner, \magnuss is taking
a random walk on the graph $\grastd{h}{f}$.  Having now phrased the
question in a precise manner, we need only give two definitions to
state the answer.

\begin{dfn}
  The \emph{Stirling number of the second kind}, $\sti{a}{b}$, counts
  the number of ways of partitioning an $a$-element set into $b$
  blocks (irrespective of order).
\end{dfn}
For example, $\sti{4}{2} = 7$ as we can partition $\{A,B,C,D\}$ into
two parts in seven ways:
  \begin{center}
    {\scalebox{.4}{\includegraphics{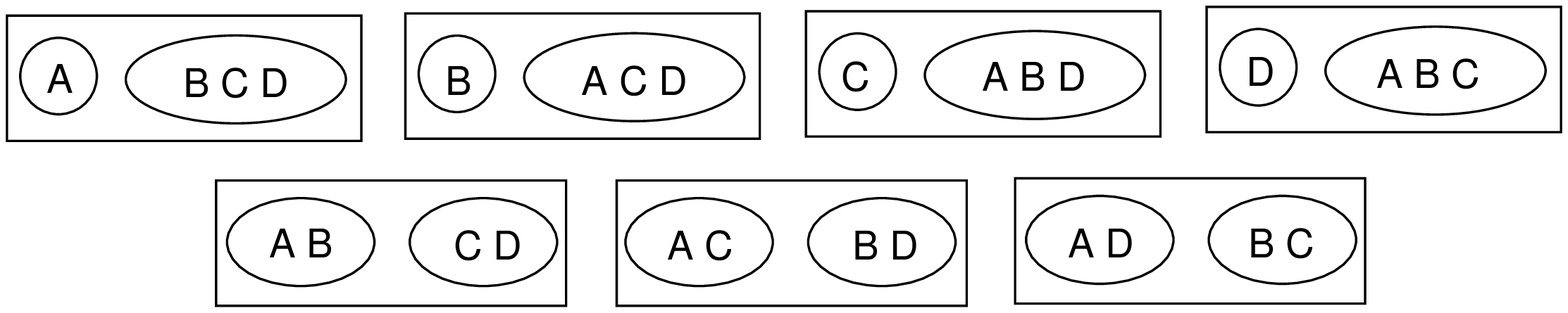}}}
  \end{center}

\begin{table}[h]
  \begin{center}
    \begin{tabular}[c]{|c|cccccc|}\hline
        & 1 & 2  & 3  & 4  & 5  & 6\\\hline
      1 & 1 &    &    &    &    &\\
      2 & 1 & 1  &    &    &    &\\
      3 & 1 & 3  & 1  &    &    &\\
      4 & 1 & 7  & 6  & 1  &    &\\
      5 & 1 & 15 & 25 & 10 & 1  &\\
      6 & 1 & 31 & 90 & 65 & 15 & 1\\\hline
    \end{tabular}
  \end{center}\label{tab:sti}
  \caption{Stirling numbers of the second kind $\sti{a}{b}$.}
\end{table}
  
Let $\stv \in\stn$.  For each integer $t$ with $1\leq t\leq h$, define
\begin{equation}
  \phi_t(\stv) = 
  \begin{cases}
    |\{t < j \leq h: \stvr_j = \emp\}|, & \text{ if } \stvr_t = \cat,\\
    0, & \text{ if } \stvr_t = \emp.
  \end{cases}
\end{equation}

\begin{dfn}
  The \emph{weight} of a state $\stv$, $\wt(\stv)$, is given by the
  formula
\begin{equation}
  \wt(\stv) = \prod_{t=1}^h (1 + \phi_t(\stv)).
\end{equation}
\end{dfn}
The weight (our terminology) arises when counting permutations with
restricted positions.  See, e.g., \cite[2.4]{S4}.


We are finally in a position to answer our Random Question.
\begin{thm}\label{thm:main}
  If \magnuss is juggling $b$ balls and is able to throw them to heights
  at most $h = b+f$, then over the long term he spends
  \begin{equation}\label{eq:fractime}
    \frac{\wt(\stv)}{\sti{h+1}{f+1}}
  \end{equation}
  of his time in the state $\stv\in\stnk$.
\end{thm}

\begin{exmp}
  Below we present the vector $\balph$ whose $\stv$-th component
  specifies the fraction of the time \magnuss spends in each state
  $\stv$ of \figref{fig:g35} while juggling 3 balls at a
  maximum throw height of 5.
\begin{equation*}
  \begin{split}
    \balph &= \frac{1}{\sti{6}{3}}{\Large(}
    \wt(\xxxoo\,),\wt(\xxoxo\,),\wt(\xxoox\,),\\
    &
    \phantom{{=} \frac{1}{\sti{6}{3}}\large(}
    \wt(\xoxxo\,),\wt(\xoxox\,),\wt(\xooxx\,),\wt(\oxxxo\,),\\
    & 
    \phantom{{=} \frac{1}{\sti{6}{3}}\large(}
    \wt(\oxxox\,), \wt(\oxoxx\,),\wt(\ooxxx\,){\Large)}\\
    &= \frac{1}{90}\left(27,18,9,12,6,3,8,4,2,1\right).
  \end{split}
\end{equation*}
From $\balph$, we see that for not quite one third of the time,
\magnuss finds himself in the state corresponding to the usual three
ball ``cascade'' juggling pattern.
\end{exmp}

\begin{remark}
  For other contexts in which the Stirling numbers of the second kind
  arise in conjunction with siteswap notation, see \cite{ER,kamstra}.
  A related class of numbers, the Eulerian numbers, also arise in the
  enumeration of siteswap patterns (see \cite{BEGW}).
\end{remark}

\section{Proof of \thmref{thm:main}}
\label{sec:main}

The proof of \thmref{thm:main} is not hard once we model the act of
juggling randomly as a \emph{Markov chain}: A discrete-time random
process in which the transition probabilities depend only on the
current state.  To describe a Markov chain, we need to know the
possible states $\Omega$ and the possible transitions among them.  The
latter we encode in a transition matrix $P$.  We will denote this pair
by $\mc{\Omega}{P}$.  The $(i,j)$-th entry of the matrix $P$ is to be
interpreted as the probability of transitioning from the $i$-th state
to the $j$-th state in one step.

In describing \magnuss as taking a random walk on the states $\stnk$
with throws chosen uniformly randomly, we have already described one
Markov chain: $\mc{\stnk}{P}$ where $P
=(p_{\stv,\stw})_{\stv,\stw\in\stnk}$ with
\begin{equation}\label{eq:tmat}
  p_{\stv,\stw} = 
  \begin{cases}
    1, & \text{ if } (\stv,\stw) \in \edges(\grastd{h}{f})\text{ and } \snuo = \emp,\\
    \frac{1}{f+1}, & \text{ if } (\stv,\stw) \in
    \edges(\grastd{h}{f})\text{ and } \snuo = \cat,\\
    0, & \text{ else}.
  \end{cases}
\end{equation}

Suppose we have a Markov chain with states $\{\stw^1,\dotcoms
\stw^r\}$ and transition matrix $P = (p_{i,j})$ such that the $l$-th
power $P^l = (p_{i,j}^{(l)})$ of $P$ satisfies $p_{i,j}^{(l)} > 0$ for
$l \gg 0$.  This condition on the matrix $P$ means that, no matter how
long \magnuss has been juggling, he still has the opportunity to visit
every state of $\stnk$.  $\mc{\stnk}{P}$ satisfies this condition.  A
\emph{probability vector} is a vector for which the sum of the entries
is 1.  The standard theorem of Markov chains we shall use is the
following:
\begin{thm}(\/\cite[4.1.4]{kemeny})\label{thm:stdmar}
  Given a Markov chain as above:
  \begin{enumerate}
  \item There exists a matrix $A$ with $\lim_{l\rightarrow \infty} P^l
    = A$.
  \item Each row of $A$ is the same vector $\balph =
    (\alpha_1,\dotcoms \alpha_r)$.
  \item Each $\alpha_i > 0$.
  \item The vector $\balph$ is the unique probability vector
    satisfying $\balph P = \balph$.
  \item For any probability vector $\boldsymbol{\pi}$,
    $\boldsymbol{\pi} P^l \rightarrow \balph$ as $l\rightarrow
    \infty$.
  \end{enumerate}
\end{thm}
The number $\alpha_i$ can be interpreted as the fraction of time the
Markov chain is expected to spend in the $i$-th state.  Hence,
\thmref{thm:stdmar} tells us how to answer our question: To find the
fraction of time that \magnuss is in state $\stw$, we need only
calculate the $\stw$-th entry of the (normalized) left eigenvector for
$P$ with eigenvalue 1!

\begin{exmp}
  We illustrate these ideas using the graph $\grastd{4}{1}$ given in
  \figref{fig:g34}.  The transition matrix $P$ simply consists
  of a $\mf{1}{2}$ for every edge except the unique edge leaving
  $\oxxx$, which must be taken.  Notice that the sum of the entries in
  a row is 1.  Below we list $P$ and two of its powers.
  \setlength{\fboxsep}{-2pt}
\begin{equation*}
P = 
\begin{matrix}
        & \xxxo & \xxox & \xoxx & \oxxx\\
  \xxxo & \mf{1}{2} & \mf{1}{2} & & \\
  \xxox & \mf{1}{2} & & \mf{1}{2} & \\
  \xoxx & \mf{1}{2} & & & \mf{1}{2} \\
  \oxxx & 1 & & &
\end{matrix}
\end{equation*}

\begin{equation*}
P^2 = 
\begin{matrix}
        & \xxxo & \xxox & \xoxx & \oxxx\\
  \xxxo & \mf{1}{2} & \mf{1}{4} & \mf{1}{4} & \\
  \xxox & \mf{1}{2} & \mf{1}{4} & & \mf{1}{4} \\
  \xoxx & \mf{3}{4} & \mf{1}{4} & & \\
  \oxxx & \mf{1}{2} & \mf{1}{2} & & \\
\end{matrix}
\end{equation*}

\begin{equation*}
P^5 = 
\begin{matrix}
        & \xxxo & \xxox & \xoxx & \oxxx\\
  \xxxo & \mf{17}{32} & \mf{9}{32} & \mf{1}{8} & \mf{1}{16} \\
  \xxox & \mf{17}{32} & \mf{1}{4} & \mf{5}{32} & \mf{1}{16} \\
  \xoxx & \mf{17}{32} & \mf{1}{4} & \mf{1}{8} & \mf{3}{32} \\
  \oxxx & \mf{9}{16} & \mf{1}{4} & \mf{1}{8} & \mf{1}{16} \\
\end{matrix}
\end{equation*}

Already for $P^5$ we see that the rows appear to be converging to the
same vector.  And, indeed, by \thmref{thm:stdmar} this vector is:
\begin{align*}
  \balph &= \frac{1}{\sti{5}{2}}\left(\wt(\xxxo\,),\wt(\xxox\,),
  \wt(\xoxx\,),\wt(\oxxx\,)\right)\\
                      &= \frac{1}{15}\left(8,4,2,1\right).
\end{align*}
Just over half of \magnus' time will be spent in the state $\xxxo$\,.
In fact, this geometric distribution is what one would expect by
examination of the graph $\grastd{4}{1}$.
\end{exmp}

One can certainly prove \thmref{thm:main} by guessing the vector
$\balph$ and then proving that it satisfies the conditions of
\thmref{thm:stdmar}.  However, as is often the case in mathematics,
keeping track of additional information enables us to write a simpler
and more illuminating proof.  So, we will construct an augmented
Markov chain that is easy to analyze and from which we can obtain our
original Markov chain by grouping states.

We will augment our states by specifying \emph{which} balls are
landing when.  To this end, define a \textit{Throwing/Landing-state}
(``TL-state'') to be an $h$-tuple $ \hsnu = \text{\fbox{\kern-3pt
    \begin{tabular}{ZZZZ} $\hnu_1$ & $\hnu_2$ & $\cdots$ & $\hnu_h$
  \end{tabular}}} \in \{\emp,1,2,\dotcoms h\}^h$ such that for each 
$1\leq j,k,t\leq h$,
\begin{equation}
\hnu_t \geq t \text{ or } \hnu_t = \emp\label{eq:cond1}
\end{equation}
and
\begin{equation}
\text{ if } j \neq k, \text{ then } \hnu_j-j \neq \hnu_k - k.\label{eq:cond2}
\end{equation}
The first condition ensures that any ball currently slated to land
has, in fact, already been thrown.  The second condition ensures that
\magnuss did not have to throw two balls at once in order to get into
his current TL-state.

As is the case for states, we interpret $\hnu_t = \emp$ to mean that
no ball lands $t$ seconds in the future.  We interpret $\hnu_t = i$ as
meaning that the ball thrown $i-t$ seconds in the past will land $t$
seconds in the future.  (Alternatively, the ball landing $t$ seconds
in the future will have been in the air $i$ seconds by the time it
lands.)  A TL-state completely captures the current status of \magnus.
Let $\stnb$ denote the set of all TL-states of length $h$ and $\stnkb$
the subset of TL-states with $f$ $\emp$'s.  The throwing operators
$\adv{j}$ are extended in the obvious manner to act on TL-states.

Of course, each TL-state determines some landing state simply by
forgetting how long each ball has been in the air.  So, define a map
$\pi_h: \stnb \longrightarrow \stn$ by $\pi_h(\hsnu) =
\text{\fbox{\kern-3pt \begin{tabular}{ZZZZ} $\snuo$ & $\snut$ &
      $\cdots$ & $\snuh$\end{tabular}}}$ where
\begin{equation}
  \nu_t = 
  \begin{cases}
    \emp, & \text{ if } \hnu_t = \emp,\\
    \cat, & \text{ if } \hnu_t \in \{1,2,\dotcoms h\}.
  \end{cases}
\end{equation}

In addition to $\mc{\stnk}{P}$, we now have another natural Markov
chain to consider: $\mc{\stnkb}{\widehat{P}}$ where $\widehat{P} =
(\widehat{p}_{\hsnu,\hsw})_{\hsnu,\hsw\in\stnkb}$ with
\begin{equation}\label{eq:tmat2}
  \widehat{p}_{\hsnu,\hsw} = 
  \begin{cases}
    1, & \text{ if } \hsw = \adv{0}(\hsnu) \text{ and } \hnu_1 = \emp,\\
    \frac{1}{f+1}, & \text{ if } \hsw = \adv{j}(\hsnu)
    \text{ for some } j \text{ with } 1\leq j \leq h,\ \hnu_1 = \cat \text{ and }
    \hnu_{j+1} = \emp,\\
    0, & \text{ else}.
  \end{cases}
\end{equation}

Our approach to proving \thmref{thm:main} will be to analyze this
latter Markov chain.  The first step will be to find the vector
$\boldsymbol{\beta}$ of steady-state probabilities for
$\mc{\stnkb}{\widehat{P}}$.  The second step will be to show that
``lumping'' states of this chain recovers our original chain
$\mc{\stnk}{P}$.  By analyzing this lumping process carefully, it will
be clear how to obtain the result of \thmref{thm:main} from our
knowledge of the vector $\boldsymbol{\beta}$.

\begin{lem}
  The vector of all 1's (of length $|\stnkb|$) is a left eigenvector
  for $\widehat{P}$.
\end{lem}
\begin{proof}
  To prove this, we need simply show that the transition matrix
  $\widehat{P}$ is doubly stochastic (i.e., that the rows \emph{and}
  columns sum to 1).
  
  There are two cases to consider.  In the first, either $\hnu_t =
  \emp$ or $\hnu_t > t$ for every $t$ with $1\leq t \leq h$.  In this
  case, \magnuss waited during the most recent beat.  That means that
  the only precursor of $\hsnu$ is $\inhnuh\,$.  On the other hand, if
  $\hnu_t = t$ for some $t$, then he just made a height $t$ throw.
  Since \magnuss throws at most one ball at any given time, there is
  at most one such $t$.  By the discussion in the proof of
  \lemref{lem:lump} given below, the previous arc of this ball could
  have been any of $f+1$ different heights.  Each of these different
  heights leads to a different precursor.  Since the out-degree of any
  state $\hsnu$ with $\hnu_1 \neq \emp$ is $f+1$, each of these edges
  has weight $\frac{1}{f+1}$.  This proves that the columns of
  $\widehat{P}$ also sum to $1$.
\end{proof}

In order to fully describe the vector $\boldsymbol{\beta}$ for
$\mc{\stnkb}{\widehat{P}}$, we need to calculate $|\stnkb|$.

\begin{lem}\label{lem:sti}
  $|\stnkb| = \sti{h+1}{f+1}$.
\end{lem}
\begin{proof}
  We proceed by constructing a bijection between elements of $\stnkb$
  and partitions with $f+1$ blocks of the set $[h+1] := \{1,2,\dotcoms
  h+1\}$.  An example of the correspondence we construct is
  illustrated in \figref{fig:partition}.
  
  We start by associating a graph $\Gamma_{\hsnu}$ with vertex set
  $[h+1]$ to each $\hsnu\in\stnkb$.  For $1\leq t < i \leq h+1$, we
  set $\{t,i\}$ to be an edge of $\Gamma_{\hsnu}$ if and only if there
  was a ball thrown $h+1-i$ seconds in the past that will be landing
  $t$ seconds in the future.  In terms of $\hsnu$, this means simply
  that $\hnu_i = h+1+t-i$.
  
  Let $\lambda(\hsnu)$ be the set partition of $[h+1]$ corresponding
  to the connected components of $\Gamma_{\hsnu}$.  Our goal is to
  show that the map $\hsnu \mapsto \lambda(\hsnu)$ is a bijection onto
  partitions of $[h+1]$ with $f+1$ blocks.
  
  To begin, we examine the connected components of $\Gamma_{\hsnu}$.
  Suppose that $\{a,b\},\{a,c\}$ are distinct edges of
  $\Gamma_{\hsnu}$.  As they are both incident to the vertex $a$, all
  of $a$, $b$ and $c$ must be in the same connected component.  It
  cannot be true that $a < b,c$ as then \magnuss would be planning to
  catch two balls $a$ seconds from now.  Nor can it be true that $b,c
  < a$ as this would imply that $a$ seconds ago \magnuss threw two
  balls at once.  We conclude that either $b < a < c$ or $c < a < b$.
  This implies that $\Gamma_{\hsnu}$ is a disjoint union of chains.
  
  The number of edges in $\Gamma_{\hsnu}$ is $h-f$ since there is an
  edge for each $j$ with $\hnu_j\neq \emp$.  So, the number of
  connected components in $\Gamma_{\hsnu}$ is the number of vertices
  ($h+1$) decreased by the number of edges ($h-f$).  This yields $f+1$
  connected components.  Hence, $\lambda(\hsnu)$ is a partition of $[h+1]$
  into $f+1$ blocks as desired.

  \mymslfig{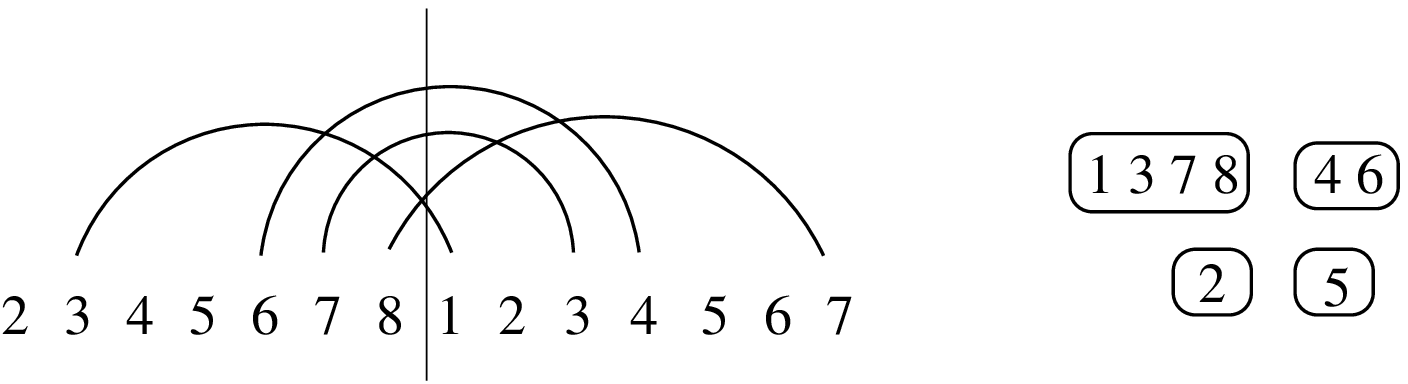}{On the left we show the edges of $\Gamma_{\hsnu}$
  for the TL-state $\hsnu = $\fbox{$\,6\smallsmile 4 6\smallsmile
    \smallsmile 7\,$}.  On the right is the corresponding partition
  $\lambda(\hsnu)$.}
  
  Now suppose that $\lambda$ is a (set) partition of $[h+1]$.  We will
  associate a TL-state 
  $\stw(\lambda) = \text{\fbox{\kern-3pt \begin{tabular}{YYY} 
        $\stwr(\lambda)_1$ & $\cdots$ & $\stwr(\lambda)_h$
      \end{tabular}}}$
  to $\lambda$.  Consider a block $\{\alpha_1,\alpha_2,\dotcoms
  \alpha_m\}$ (indexed such that $\alpha_i < \alpha_j$ if $i < j$) of
  $\lambda$.  If $m = 1$, then set $\stwr(\lambda)_{\alpha_1} = \emp$.
  Otherwise, set $\stwr(\lambda)_{\alpha_l} =
  h+1+\alpha_l-\alpha_{l+1}$ for each $l$ with $1\leq l < m$.  Since
  $1\leq \alpha_l < \alpha_{l+1} \leq h+1$ for each $l$,
  $\stw(\lambda)$ is certainly an $h$-tuple in $\{\emp,1,2,\dotcoms
  h\}$.  That \eqref{eq:cond1} is satisfied follows from the
  inequalities $\alpha_l < \alpha_{l+1} \leq h+1$.  That
  $\stw(\lambda)$ satisfies \eqref{eq:cond2} follows from the
  definition of $\stwr(\lambda)_{\alpha_l}$ along with the fact that
  our blocks are disjoint.  Hence, $\stw(\lambda)$ is, in fact, a
  TL-state.  Furthermore, it can be checked that 
  $\lambda \mapsto \stw(\lambda)$ is the inverse of the map 
  $\hsnu\mapsto \lambda(\hsnu)$.  Therefore, every partition of $[h+1]$ into $f+1$
  blocks is obtained as a $\lambda(\hsnu)$ for a unique
  $\hsnu\in\stnkb$.  This completes the proof.
\end{proof}

\begin{cor}\label{cor:uni}
  In the process $\mc{\stnkb}{\widehat{P}}$, the fraction of time
  \magnuss finds himself in state $\hsnu$ is
  $\frac{1}{\sti{h+1}{f+1}}$.
\end{cor}

We have shown that, when juggling randomly, \magnuss is as likely to
find himself in one TL-state as another.  We will now consider a new
random process running in parallel with $\mc{\stnkb}{\widehat{P}}$.
For $\stv\in\stnk$, define $[\stv] = \pi_h^{-1}(\stv)\in\stnkb$.  When
$\mc{\stnkb}{\widehat{P}}$ is in state $\hsnu$, we define our new
process to be in the state $[\pi_h(\hsnu)]$.  It follows that its
states are $\{[\stv]\}_{\stv\in\stnk}$.  This is a \emph{lumped}
process derived from $\mc{\stnkb}{\widehat{P}}$ by grouping together
certain states.


Certainly \magnuss wanders randomly among the states of this new
process.  We do not yet know, however, that this new process is a
Markov chain.  In particular, it is not clear that the probability of
transitioning from $[\stv]$ to $[\stw]$ is independent of how long
\magnuss has been in $[\stv]$.  If it \emph{were} dependent, then the
transition probabilities would depend on previous states as well as
the current state.  Or, stated another way, the transition probability
from $[\stv]$ to $[\stw]$ would depend on \emph{which} state $\hsmu$
of $[\stv]$ \magnuss is currently in.

For each $\hsmu\in[\stv]$, let $\widehat{p}_{\hsmu,[\stw]}$ denote the
probability of $\mc{\stnkb,\widehat{P}}$ transitioning from $\hsmu$ to
some element of $[\stw]$.  Our above discussion suggests (see
\cite[6.3.2]{kemeny} for a proof) that our new process will be a
Markov chain if the probabilities $p_{\hsmu,[\stw]}$ are equal for
each $\hsmu\in[\stv]$.  Furthermore, the transition probabilities
$r_{[\stv],[\stw]}$ of our new process will be these common values
$\widehat{p}_{\hsmu,[\stw]}$.

\begin{lem}
  Fix $\stv,\stw\in\stnk$.  With the notation above,
  $\widehat{p}_{\hsmu,[\stw]} = p_{\stv,\stw}$ for each
  $\hsmu\in[\stv]$.
\end{lem}
\begin{proof}
  If $p_{\stv,\stw} = 0$, then certainly all of the
  $\widehat{p}_{\hsmu,[\stw]}$ are $0$ too.  Also, if $p_{\stv,\stw} =
  1$, then \magnuss must wait after any state $\hsmu\in[\stv]$; hence
  $\widehat{p}_{\hsmu,[\stw]} = 1$.  So the only case to consider is
  when $p_{\stv,\stw} = \frac{1}{f+1}$.  In this case, $\stw =
  \adv{j}(\stv)$ for some $j$ with $1\leq j\leq h$ and $\nu_1 = \cat$.
  In particular, there is some $j$ with $1\leq j\leq h$ for which
  $\nu_{j+1} = \emp \text{ and } \omega_j = \cat$.  Suppose that
  \magnuss transitions from $\hsmu$ to a particular $\hsw\in[\stw]$.
  The previous equalities imply that $\hmu_{j+1} = \emp$ and
  $\widehat{\omega}_j = j$.  It follows that $\hsw$ is obtained
  from $\hsmu$ by making a height $j$ throw.  This will be one of
  $f+1$ possible throws for \magnuss as $\hmu_1 \neq \emp$; hence he
  will take this choice with probability $\frac{1}{f+1}$.  As we have
  considered all cases, this completes the proof of the lemma.
\end{proof}

Define a matrix $R = (r_{[\stv],[\stw]})_{\stv,\stw\in\stnk}$ by
setting $r_{[\stv],[\stw]} = p_{\stv,\stw}$.  By the above lemma, we
have obtained a new stochastic process
$\mc{\{[\stv]\}_{\stv\in\stnk}}{R}$ that can be identified in the
obvious manner with $\mc{\stnk}{P}$.  It follows from \corref{cor:uni}
that the fraction of the time \magnuss finds himself in a state $\stv$
with respect to the process $\mc{\stnk}{P}$ is
$\frac{|[\stv]|}{\sti{h+1}{f+1}}$.  Hence, \thmref{thm:main} follows
from \lemref{lem:lump}.

\mymslfig{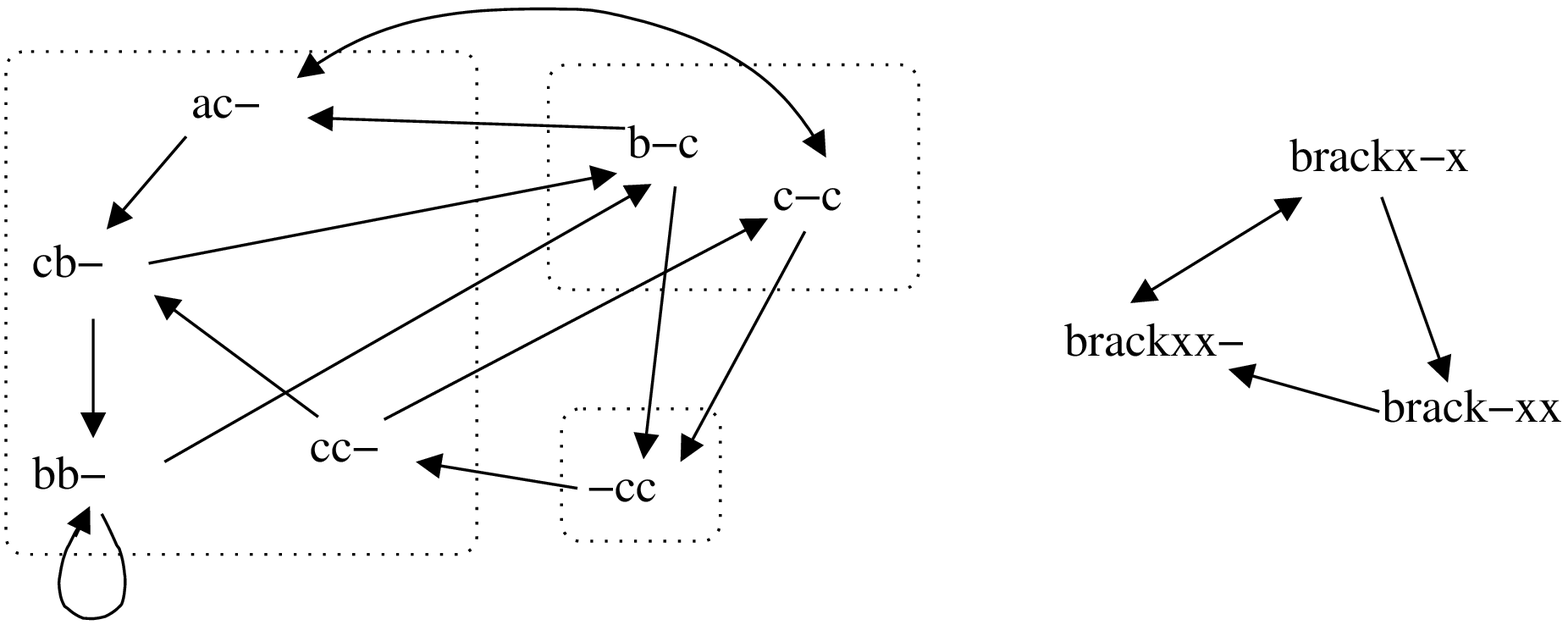}{The Markov chain
  $\mc{\widehat{St}_{3,1}}{\widehat{P}}$ along with the lumped
  process.}

\begin{lem}\label{lem:lump}
  \begin{equation}\label{eq:invcard}
    |[\stv]| = \wt(\stv) = \prod_{t=1}^h (1+\phi_t(\stv)).
  \end{equation}
\end{lem}
\begin{proof}
  Let $\hsnu\in\stnkb$ and $\stv\in\stnk$ such that $\pi_h(\hsnu) =
  \stv$.  Suppose $\nu_t = \cat$ and that we know the value of
  $\hnu_i$ for each $i > t$.  How do these values constrain the
  possibilities for $\hnu_t$?  To start, note that the ball landing
  $t$ seconds in the future was thrown at most $h-t$ seconds in the
  past.  As \magnuss must have already thrown this ball, there are at
  most $h-t+1$ possibilities for $\hnu_t$.  In addition, any ball
  landing \emph{more} than $t$ seconds in the future must have been
  thrown \emph{fewer} than $h-t$ seconds in the past as $\stv\in\stn$.
  Combining these two observations, we see that the number of
  possibilities for $\hnu_t$ is reduced precisely by the number of
  balls landing after it.  In particular, the possible values for
  $\hnu_t$ depend only on the $\stvr_i$ for $i > t$.
  
  In conclusion, the number of possible values for $\hnu_t$ when
  $\stvr_t = \cat$ is $$h-t+1 - |\{i > t: \stvr_i = \cat\}|,$$
  which
  is just $(1+\phi_t(\stv))$.  Taking the product over all $t$ yields
  the desired answer.
\end{proof}

This completes the proof of \thmref{thm:main}.

\begin{remark}
  Consideration of the powers of the transition matrix $P$ might lead
  to new asymptotic formulae for the Stirling numbers of the second
  kind: Consider $\grastd{h}{f}$ and the associated transition matrix
  $P$ with the last row indexed by the state $\stw = \odoxdx\in\stnk$.
  We have already seen that the rows of $P^l$ converge to a
  probability vector $\balph$ with $\alpha_{\stw} = 1/\sti{h+1}{f+1}$.
  Hence, for each state $\stu\in\stnk$, we get a sequence
  \begin{equation*}
    (p_{\stu,\stw},p_{\stu,\stw}^{(2)},\ldots)
  \end{equation*}
  that converges to $1/\sti{h+1}{f+1}$.  
\end{remark}

\section{Other models}
\label{sec:models}

By now (unless you are a very fast reader), \magnuss is getting tired.
He is going to drop occasionally.  Certainly any realistic juggling
model should accommodate transitions to states with fewer balls.  Of
course, if we allow \magnuss to juggle indefinitely, and he drops
occasionally, he will eventually run out of balls.  To keep things
interesting, we will give his assistant, Sphagnum, the option of
inserting a ball into \magnus' pattern whenever \magnuss has a wait.

We will explore two such models in this section.  In both models,
\magnuss can simply fail to catch a ball.  When this happens, \magnuss
transitions from a state in $\stnk$ to one in $\stnkp$.  The
difference between the two models lies in what is considered to be a
legal throw.  In the first model, \magnuss (or his assistant) will
only throw to heights that would not lead to two balls landing at the
same time.  In the second model, however, we allow this eventuality.
Of course, as \magnuss is well-known to have small hands, he cannot
catch two at once, so he will have to ignore one of the balls.  Such a
throw will, in terms of transitions, be treated no differently than a
drop.

These changes are easily incorporated into the framework we have
developed.  Now we define the state graphs corresponding to the
aforementioned \emph{add-drop model} (well known to college
students) and \emph{annihilation model}.

\begin{enumerate}
\item \textit{Add-drop juggling} generalizes standard juggling by
  always allowing throws of height 0 and removing the restriction on
  $\nu_1 = \cat$ for throws of positive height.  (The throws of height
  $0$ correspond to \magnuss dropping or waiting; those of positive
  height when $\nu_1 = \emp$ correspond to Sphagnum helping.) So, our
  graph $\graddd{h}$ has
  \begin{enumerate}
  \item Vertices indexed by $\stn$.
  \item An edge from $\stv$ to $\stw$ whenever 
     \begin{enumerate}
     \item $\stw = \adv{j}(\stv)$ for some $j$, $1 \leq j \leq h$ such
       that $\nu_{j+1} = \emp$; or
     \item $\stw = \adv{0}(\stv)$.
     \end{enumerate}
  \end{enumerate}
  \figref{fig:ad3} illustrates $\graddd{3}$.
  \mymslfig{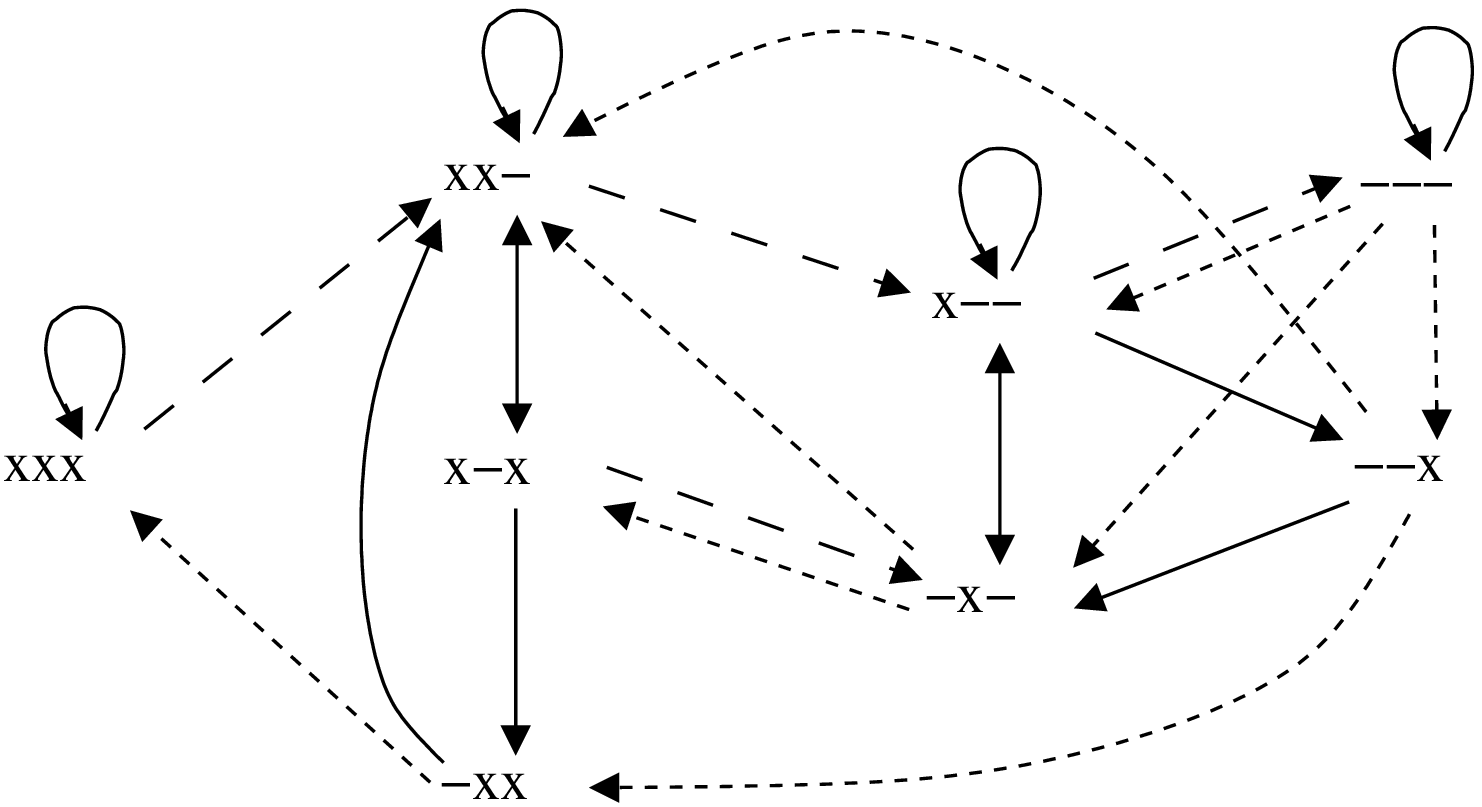}{Solid arrows are regular throws; dashed are
    drops; dotted are insertions by Sphagnum.}
  
\item In \textit{Annihilation juggling}, we take the same state graph
  $\grann{h} = \graddd{h}$ as in the previous case.  The difference will
  be in how we assign probabilities to the edges.
\end{enumerate}

With our state graphs now defined, we can define the corresponding
transition matrices.  For $\stv\in\stnk$, set 
$f'_{\stv} = |\{t: 2\leq t\leq h \text{ and } \nu_t = \emp\}|$.
\begin{enumerate}
\item For $\graddd{h}$, set $Q = (q_{\stv,\stw})_{\stv,\stw\in\stn}$ with
  \begin{equation*}
    q_{\stv,\stw} = 
    \begin{cases}
      \frac{1}{f'_{\stv} + 2}, & \text{ if } (\stv,\stw) \in \edges(\graddd{h}),\\
      0, & \text{ else}.
    \end{cases}
  \end{equation*}
  The denominator above comes from the fact that \magnuss can either
  drop/wait or throw to any $\emp$ occurring at $\nu_2$ or later.  So
  here, again, \magnuss chooses uniformly from the available edges.

\item For $\grann{h}$, set $R = (r_{\stv,\stw})_{\stv,\stw\in\stn}$ with
  \begin{equation*}
    r_{\stv,\stw} = 
    \begin{cases}
      \frac{h-f'_{\stv}}{h+1}, & \text{ if } (\stv,\stw) \in \edges(\grann{h})
      \text{ and } \adv{0}(\stv) = \stw,\\
      \frac{1}{h+1}, & \text{ if } (\stv,\stw) \in \edges(\grann{h})
      \text{ and } \adv{0}(\stv) \neq \stw,\\
      0, & \text{ else}.
    \end{cases}
  \end{equation*}
  Here, you should envision \magnuss (possibly with the help of
  Sphagnum) throwing to \emph{any} height ($\leq h$) at random,
  regardless of whether or not a ball has just landed.  If there is
  already a ball landing at the time he is throwing to, the latter
  ball gets annihilated.
\end{enumerate}

In order to describe the steady-state probabilities for the add-drop
model, we must introduce the Bell numbers.  The $h$-th Bell number,
$B_h$, is defined by $ B_h = \sum_{i=0}^h \sti{h}{i}$.  For these two
models, the analogue of \thmref{thm:main} is
\begin{thm}\label{thm:other}
  Let $\stv\in \stnk$.
  \begin{enumerate}
  \item In the add-drop model, \magnuss spends 
    \begin{equation*}
      \frac{\wt(\stv)}{\bell{h+1}}
    \end{equation*}
    of his time in state $\stv$.
  \item In the annihilation model, \magnuss spends
    \begin{equation*}
      \frac{\frac{h!}{(h-f)!}\wt(\stv)}{(h+1)^h}
    \end{equation*}
    of his time in state $\stv$.
  \end{enumerate}
\end{thm}

As the denominators for our probabilities in the standard and add-drop
model have natural combinatorial interpretations, the reader may be
wondering if such an interpretation exists for the $(h+1)^h$ occurring
in the annihilation model.  Indeed, $(h+1)^h$ counts the number of
labeled rooted trees with $h$ nodes (see \cite[5.3.2]{ECII} for details).

\begin{exmp}
Below we give the transition matrix for $\graddd{3}$:
\setlength{\fboxsep}{-2pt}
\begin{equation*}
  \begin{matrix}
    &
\text{\fbox{\kern-2pt\begin{tabular}{YYY} $\cat$ & $\cat$ & $\cat$\end{tabular}}}&
\text{\fbox{\kern-2pt\begin{tabular}{YYY} $\cat$ & $\cat$ & $\emp$\end{tabular}}}&
\text{\fbox{\kern-2pt\begin{tabular}{YYY} $\cat$ & $\emp$ & $\cat$\end{tabular}}}&
\text{\fbox{\kern-2pt\begin{tabular}{YYY} $\emp$ & $\cat$ & $\cat$\end{tabular}}}&
\text{\fbox{\kern-2pt\begin{tabular}{YYY} $\cat$ & $\emp$ & $\emp$\end{tabular}}}&
\text{\fbox{\kern-2pt\begin{tabular}{YYY} $\emp$ & $\cat$ & $\emp$\end{tabular}}}&
\text{\fbox{\kern-2pt\begin{tabular}{YYY} $\emp$ & $\emp$ & $\cat$\end{tabular}}}&
\text{\fbox{\kern-2pt\begin{tabular}{YYY} $\emp$ & $\emp$ & $\emp$\end{tabular}}}\\
\text{\fbox{\kern-2pt\begin{tabular}{YYY} $\cat$ & $\cat$ & $\cat$\end{tabular}}}&
    \mf{1}{2} & \mf{1}{2} & & & & & & \\
\text{\fbox{\kern-2pt\begin{tabular}{YYY} $\cat$ & $\cat$ & $\emp$\end{tabular}}}&
    & \mf{1}{3} & \mf{1}{3} & & \mf{1}{3} & & & \\
\text{\fbox{\kern-2pt\begin{tabular}{YYY} $\cat$ & $\emp$ & $\cat$\end{tabular}}}&
    & \mf{1}{3} & & \mf{1}{3} & & \mf{1}{3} & & \\
\text{\fbox{\kern-2pt\begin{tabular}{YYY} $\emp$ & $\cat$ & $\cat$\end{tabular}}}&
    \mf{1}{2} & \mf{1}{2} & & & & & & \\
\text{\fbox{\kern-2pt\begin{tabular}{YYY} $\cat$ & $\emp$ & $\emp$\end{tabular}}}&
    & & & & \mf{1}{4} & \mf{1}{4} & \mf{1}{4} & \mf{1}{4} \\
\text{\fbox{\kern-2pt\begin{tabular}{YYY} $\emp$ & $\cat$ & $\emp$\end{tabular}}}&
    & \mf{1}{3} & \mf{1}{3} & & \mf{1}{3} & & & \\
\text{\fbox{\kern-2pt\begin{tabular}{YYY} $\emp$ & $\emp$ & $\cat$\end{tabular}}}&
    & \mf{1}{3} & & \mf{1}{3} & & \mf{1}{3} & & \\
\text{\fbox{\kern-2pt\begin{tabular}{YYY} $\emp$ & $\emp$ & $\emp$\end{tabular}}}&
    & & & & \mf{1}{4} & \mf{1}{4} & \mf{1}{4} & \mf{1}{4} \\
  \end{matrix}
\end{equation*}
\end{exmp}

From \thmref{thm:other}.1, we obtain the probability vector 
\begin{equation*}
  \boldsymbol{\balph} = \frac{1}{15}\left(1,4,2,1,3,2,1,1\right).
\end{equation*}
This gives the fraction of the time \magnuss spends in each of the
eight possible states (the ordering of $\balph$ is the same as that of
the transition matrix).

As the steps in the proof of \thmref{thm:other} are analogous to those
found in the proof of \thmref{thm:main}, we do not detail them here.
As a consolation to the reader, though, we mention two other families
of juggling models which, although perhaps less physically realistic
than the ones we have described, do lead to interesting answers to our
Random Question.

\begin{enumerate}
\item \emph{Multiplex juggling}: There is no real reason (assuming you
  find a juggler with bigger hands) to disallow the catching and
  throwing of more than one ball at a time.  The standard and add-drop
  models can both be reinterpreted with this condition relaxed.
\item \emph{Infinite juggling}: Again, given a juggler in better shape
  than \magnus, there is no reason to place a maximum on legal throw
  heights.  In the standard juggling model, placing a probability of
  $(p-1)/p^j$ for the throw to the $j$-th next available height yields
  tractable calculations for any $p > 1$.  
  
  Rather than choosing geometric weights for the throws, one could
  choose probabilities from any convergent infinite series.  For
  example, why not choose $6/(\pi^2 j^2)$ for the $j$-th available
  throw height?  Or even $90/(\pi^4 j^4)$?  The sky's the limit.
\end{enumerate}

\section{Acknowledgements}
I would like to thank Michael Schneider for helpful advice on Markov
chains and Joe Buhler for help in sorting out who should be credited
with the idea of siteswap notation.

\bibliography{gen}
\end{document}

%% file: macros.tex


\hfuzz4pt 

\setlength{\floatsep}{15pt plus 12pt}
\setlength{\textfloatsep}{\floatsep}

\setlength{\intextsep}{\floatsep}\def\ldotsplus{\mathinner{\ldotp\ldotp\ldotp\ldotp}}
\setlength{\intextsep}{\floatsep}\def\ldotscomm{\mathinner{\ldotp\ldotp\ldotp,}}
\def\fourdots{\relax\ifmmode\ldotsplus\else$\m@th \ldotsplus\,$\fi}
\def\dotcoms{\relax\ifmmode\ldotscomm\else$\m@th \ldotsplus\,$\fi}

%
%
\setlength{\unitlength}{0.06em}
\newlength{\cellsize} \setlength{\cellsize}{18\unitlength}
\newsavebox{\cell}
\sbox{\cell}{\begin{picture}(18,18)
\put(0,0){\line(1,0){18}}
\put(0,0){\line(0,1){18}}
\put(18,0){\line(0,1){18}}
\put(0,18){\line(1,0){18}}
\end{picture}}
\newcommand\cellify[1]{\def\thearg{#1}\def\nothing{}%
\ifx\thearg\nothing
\vrule width0pt height\cellsize depth0pt\else
\hbox to 0pt{\usebox{\cell} \hss}\fi%
\vbox to \cellsize{
\vss
\hbox to \cellsize{\hss$#1$\hss}
\vss}}
\newcommand\tableau[1]{\vtop{\let\\\cr
\baselineskip -16000pt \lineskiplimit 16000pt \lineskip 0pt
\ialign{&\cellify{##}\cr#1\crcr}}}
%

\theoremstyle{plain}
\newtheorem{thm}{Theorem}
\newtheorem{lem}[thm]{Lemma}
\newtheorem{cor}[thm]{Corollary}
\newtheorem{prop}[thm]{Proposition}
\newtheorem{conj}[thm]{Conjecture}

\theoremstyle{definition}
\newtheorem{dfn}[thm]{Definition}
\newtheorem{exmp}[thm]{Example}
\newtheorem{remark}[thm]{Remark}
\newtheorem{fact}[thm]{Fact}



\newcommand{\thmref}[1]{Theorem~\ref{#1}}
\newcommand{\corref}[1]{Corollary~\ref{#1}}
\newcommand{\conjref}[1]{Conjecture~\ref{#1}}
\newcommand{\propref}[1]{Proposition~\ref{#1}}
\newcommand{\propsref}[1]{Propositions~\ref{#1}}
\newcommand{\lemref}[1]{Lemma~\ref{#1}}
\newcommand{\lemsref}[1]{Lemmas~\ref{#1}}
\newcommand{\factref}[1]{Fact~\ref{#1}}
\newcommand{\rmkref}[1]{Remark~\ref{#1}}
\newcommand{\remarkref}[1]{Remark~\ref{#1}}
\newcommand{\defref}[1]{Definition~\ref{#1}}
\newcommand{\egref}[1]{Example~\ref{#1}}
\newcommand{\figref}[1]{Figure~\ref{#1}}
\newcommand{\figsref}[2]{Figures~\ref{#1},\ref{#2}}
\newcommand{\secref}[1]{Section~\ref{#1}}

\newcommand{\cwd}{cwd}
\newcommand{\rwd}{rwd}
\newcommand{\ltx}{P_x}
\newcommand{\rtx}{Q_x}
\newcommand{\ltw}{P_w}
\newcommand{\rtw}{Q_w}
\newcommand{\rtq}{Q}
\newcommand{\rtqp}{Q'}
\newcommand{\rtqpp}{Q''}
\newcommand{\bigvert}{!{\vrule width 1pt}}

\newcommand{\miniprefig}
      {\scalebox{.25}{\includegraphics{hexpics/ministr.eps}}}

\newcommand{\minimslfig}
      {\scalebox{.08}{\includegraphics{mslpics/sec3412.eps}}}

\newcommand{\miniinvfig}
      {\scalebox{.08}{\includegraphics{mslpics/sec4231.eps}}}

\newcommand{\minihexafig}
      {\scalebox{.07}{\includegraphics{mslpics/minihexa.eps}}}
\newcommand{\minihexbfig}
      {\scalebox{.07}{\includegraphics{mslpics/minihexb.eps}}}
\newcommand{\minihexcfig}
      {\scalebox{.07}{\includegraphics{mslpics/minihexc.eps}}}
\newcommand{\minihexdfig}
      {\scalebox{.07}{\includegraphics{mslpics/minihexd.eps}}}
\newcommand{\minihexefig}
      {\scalebox{.07}{\includegraphics{mslpics/minihexe.eps}}}
\newcommand{\minihexffig}
      {\scalebox{.07}{\includegraphics{mslpics/minihexf.eps}}}
\newcommand{\minihexgfig}
      {\scalebox{.07}{\includegraphics{mslpics/minihexg.eps}}}

\newcommand{\mymslfig}[2]{\begin{figure}[htbp]\begin{center}
      {\scalebox{.6}{\includegraphics{#1}}}
      \caption{#2}\label{fig:#1}
    \end{center}\end{figure}}

\newcommand{\myintromslfig}[2]{\begin{figure}[htbp]\begin{center}
      {\scalebox{.4}{\includegraphics{#1.eps}}}
      \footnotesize{\textbf{Figure 0.1.} #2}
    \end{center}\end{figure}}

\newcommand{\mysmmslfig}[2]{\begin{figure}[ht]\begin{center}
      {\scalebox{.32}{\includegraphics{#1.eps}}}
      \caption{#2}\label{fig:#1}
    \end{center}\end{figure}}

\newcommand{\myhexfig}[2]{\begin{figure}[ht]\begin{center}
    \input{hexpics/#1.pstex_t}\caption{#2}\label{fig:#1}
    \end{center}\end{figure}}

\newcommand{\mynewhexfig}[2]{\begin{figure}[ht]\begin{center}
    \input{hexpics/#1.pstex_t}\caption{#2}\label{fig:#1}
    \end{center}\end{figure}}

\newcommand{\mynewfig}[2]{\begin{figure}[ht]\begin{center}
      {\scalebox{.3}{\includegraphics{newpics/#1.eps}}}
      \caption{#2}\label{fig:#1}
    \end{center}\end{figure}}

\newcommand{\sfrac}[2]{\genfrac{\{}{\}}{0pt}{}{#1}{#2}}
\newcommand{\sumsb}[1]{\sum_{\substack{#1}}}  
\newcommand{\minsb}[1]{\substack{#1}}  
\newcommand{\pstack}[1]{P_{\minsb{#1}}}  
\newcommand{\musumr}[3]{\Theta_{(\cdot #3)}[#1,#2]}
\newcommand{\musuml}[3]{\Theta_{(#3 \cdot)}[#1,#2]}
\newcommand{\mscor}[3]{\delta_{(\cdot #3)}[#1,#2]}
\newcommand{\mslor}[3]{\omega_{(\cdot #3)}[#1,#2]}
\newcommand{\msset}[3]{\theta_{(\cdot #3)}[#1,#2]}
        
\newcommand{\mscol}[3]{\delta_{(#3 \cdot)}[#1,#2]}
\newcommand{\mslol}[3]{\omega_{(#3 \cdot)}[#1,#2]}

\newcommand{\fsk}{S_k}
\newcommand{\fsn}{S_n}
\newcommand{\fsp}[1]{S_{#1}}
\newcommand{\frg}{\mathfrak{g}}
\newcommand{\frh}{\mathfrak{h}}
\newcommand{\slgn}{\mathfrak{sl}_n}

\newcommand{\io}{{i_1}}
\newcommand{\iw}{{i_2}}
\newcommand{\inn}{{i_r}}

\newcommand{\bbA}{\mathbb{A}}
\newcommand{\bbR}{\mathbb{R}}
\newcommand{\bbC}{\mathbb{C}}
\newcommand{\bbZ}{\mathbb{Z}}
\newcommand{\bbN}{\mathbb{N}}
\newcommand{\bbQ}{\mathbb{Q}}

\newcommand{\wbi}{w^{\hat{i}}}
\newcommand{\xbi}{x^{\hat{i}}}
\newcommand{\zbi}{z^{\hat{i}}}
\newcommand{\wbio}{w^{\hat{1}}}
\newcommand{\wsbio}{ws^{\hat{1}}}
\newcommand{\xbio}{x^{\hat{1}}}
\newcommand{\zbio}{z^{\hat{1}}}

\newcommand{\tti}{t}
\newcommand{\wti}{{\widetilde{w}}}
\newcommand{\xti}{\protect{\widetilde{x}}}
\newcommand{\uti}{{\widetilde{u}}}
\newcommand{\vti}{{\widetilde{v}}}
\newcommand{\yti}{{\widetilde{y}}}
\newcommand{\zti}{{\widetilde{z}}}

\newcommand{\nts}{\negthickspace}

\newcommand{\trans}{\mathcal{T}}
\newcommand{\simref}{\mathcal{S}}
\newcommand{\calh}{\mathcal{H}}
\newcommand{\calw}{\mathcal{W}}


\newcommand{\msp}{\textsc{msp}}

\newcommand{\br}{\mathbf{a}}
\newcommand{\brp}{\mathbf{a'}}
\newcommand{\bsig}{{\boldsymbol{\sigma}}}
\newcommand{\bnu}{{\boldsymbol{\nu}}}
\newcommand{\bgam}{{\boldsymbol{\gamma}}}
\newcommand{\bdel}{{\boldsymbol{\delta}}}
\newcommand{\bal}{{\boldsymbol{\alpha}}}
\newcommand{\bbe}{{\boldsymbol{\beta}}}

\newcommand{\rds}{\operatorname{D_R}}
\newcommand{\lds}{\operatorname{D_L}}
\newcommand{\musum}{\operatorname{Musum}}
\newcommand{\rank}{\operatorname{lvl}}
\newcommand{\pt}{\operatorname{pt}}
\newcommand{\heap}{\operatorname{Heap}}
\newcommand{\lcz}{\operatorname{lcz}}
\newcommand{\rcz}{\operatorname{rcz}}
\newcommand{\mcz}{\operatorname{mcz}}
\newcommand{\ver}{\operatorname{ver}}
\newcommand{\codim}{\operatorname{codim}}
\newcommand{\bcone}{\operatorname{Cone_{\wedge}}}
\newcommand{\ucone}{\operatorname{Cone^{\vee}}}
\newcommand{\charp}{\operatorname{char}}

\newcommand{\IH}{\operatorname{IH}}
\newcommand{\ext}{\operatorname{Flush}}

\newcommand{\defeq}{:=}

\newcommand{\digx}{\operatorname{mat}(x)}
\newcommand{\digxt}{\operatorname{mat}(\widetilde{x})}
\newcommand{\digwt}{\operatorname{mat}(\widetilde{w})}
\newcommand{\augdigx}{\operatorname{mat}'(x)}
\newcommand{\digw}{\operatorname{mat}(w)}
\newcommand{\digp}[1]{\operatorname{mat}(#1)}
\newcommand{\dpw}[1]{d_{#1,w}}
\newcommand{\dpb}[2]{d_{#1,#2}}
\newcommand{\duv}{d_{u,v}}
\newcommand{\dvw}{d_{v,w}}
\newcommand{\dzw}{d_{z,w}}
\newcommand{\dxw}{d_{x,w}}
\newcommand{\dyw}{d_{y,w}}
\newcommand{\txw}{\Theta_{x,w}}
\newcommand{\dxtwt}{d_{\widetilde{x},\widetilde{w}}}

\newcommand{\rxjw}{\mathcal{R}(x_i,w)}
\newcommand{\rxjmw}{\mathcal{R}(x_{i-1},w)}
\newcommand{\rxow}{\mathcal{R}(x_0,w)}

\newcommand{\area}{\mathcal{A}}

\newcommand{\ph}{\phi}
\newcommand{\pht}{\phi_t}
\newcommand{\phs}{\phi_s}
\newcommand{\phtp}{\phi_{t'}}
\newcommand{\phtal}{\phi_{t_{\alpha_j}}}

\newcommand{\cP}{{\mathcal P}}
\newcommand{\xsing}{{X_w^{\text{sing}}}}
\newcommand{\maxsing}{\operatorname{\protect{Maxsing}}(X_\protect{w})}
\newcommand{\maxsingt}{\operatorname{Maxsing}(X_{\widetilde{w}})}
\newcommand{\schub}[1]{X_{#1}}

\newcommand{\sij}{s_{i_j}}
\newcommand{\sik}{s_{i_k}}
\newcommand{\sio}{s_{i_1}}
\newcommand{\sir}{s_{i_r}}


\newcommand{\puv}{P_{u,v}}
\newcommand{\pxw}{P_{x,w}}
\newcommand{\pww}{P_{w,w}}
\newcommand{\pyw}{P_{y,w}}
\newcommand{\pyv}{P_{y,v}}
\newcommand{\pzv}{P_{z,v}}
\newcommand{\pxws}{P_{x,ws}}
\newcommand{\pxsws}{P_{xs,ws}}
\newcommand{\pxtwt}{P_{\xti,\wti}}
\newcommand{\pzw}{P_{z,w}}
\newcommand{\pxiwi}{P_{x^{-1},w^{-1}}}
\newcommand{\psxw}{P_{sx,w}}
\newcommand{\pxsw}{P_{xs,w}}
\newcommand{\pxz}{P_{x,z}}
\newcommand{\pkl}[2]{P_{#1,#2}}


\newcommand{\delxw}{\Delta(x,w)}
\newcommand{\delyw}{\Delta(y,w)}
\newcommand{\delytw}{\Delta(yt,w)}

\newcommand{\cpw}{{\mathcal P}(\br)}
\newcommand{\cpww}{{\mathcal P}_w(\br)}
\newcommand{\cpws}{{\mathcal P}(\br s)}
\newcommand{\cpxw}{{\mathcal P}_x(\br)}
\newcommand{\cpyw}{{\mathcal P}_y(\br)}
\newcommand{\cpxsw}{{\mathcal P}_{xs}(\br)}
\newcommand{\cpxsws}{{\mathcal P}_{xs}(\br/s)}
\newcommand{\cpxws}{{\mathcal P}_x(\br/s)}
\newcommand{\cpzero}{{\mathcal P}_x^0(\br)}
\newcommand{\cpone}{{\mathcal P}_x^1(\br)}
\newcommand{\cpeps}{{\mathcal P}_x^\epsilon(\br)}

\newcommand{\cd}{{\mathcal D}}
\newcommand{\Dbr}{\Delta_\bsig}
\newcommand{\Dbrj}{\Delta_{\bsig[j]}}
\newcommand{\Dbrr}{\Delta_{\bsig[r]}}
\newcommand{\Dbrk}{\Delta_{\bsig[k]}}
\newcommand{\Dbrjm}{\Delta_{\bsig[j-1]}}

\newcommand{\apw}{P(\br)}
\newcommand{\apww}{P_w(\br)}
\newcommand{\apxw}{P_x(\br)}
\newcommand{\apxsw}{P_{xs}(\br)}
\newcommand{\apxws}{P_x(\br/s)}
\newcommand{\apxsws}{P_{xs}(\br/s)}

\newcommand{\apzw}{P_z(\br)}
\newcommand{\apzsw}{P_{zs}(\br)}
\newcommand{\apzws}{P_z(\br/s)}
\newcommand{\apzsws}{P_{zs}(\br/s)}

\newcommand{\apew}{P_e(\br)}
\newcommand{\cpew}{{\mathcal P}_e(\br)}


\newcommand{\fl}{\operatorname{fl}}
\newcommand{\ra}{\operatorname{unfl}}
\newcommand{\im}{\operatorname{Im}}
\newcommand{\slide}{\operatorname{\tau}}

\newcommand{\tpq}{t_{p,q}}
\newcommand{\tab}{t_{a,b}}
\newcommand{\talbt}{t_{\alpha,\beta}}
\newcommand{\tac}{t_{a,c}}
\newcommand{\tbc}{t_{b,c}}
\newcommand{\tcd}{t_{c,d}}
\newcommand{\tij}{t_{i,j}}
\newcommand{\tik}{t_{i,k}}
\newcommand{\tjk}{t_{j,k}}
\newcommand{\til}{t_{i,l}}
\newcommand{\talj}{t_{\alpha_j}}
\newcommand{\tgd}{t_{\gamma,\delta}}

\newcommand{\ttpq}{t_{p,q}}
\newcommand{\ttab}{t_{a,b}}
\newcommand{\ttalbt}{t_{\alpha,\beta}}
\newcommand{\ttac}{t_{a,c}}
\newcommand{\ttbc}{t_{b,c}}
\newcommand{\ttcd}{t_{c,d}}
\newcommand{\ttij}{t_{i,j}}
\newcommand{\ttik}{t_{i,k}}
\newcommand{\ttjk}{t_{j,k}}
\newcommand{\ttalj}{t_{\alpha_j}}
\newcommand{\ttgd}{t_{\gamma,\delta}}

\newcommand{\rp}[2]{\mathcal{R}({#1},{#2})}
\newcommand{\qp}[2]{\mathcal{E}_{{#1},{#2}}(x,w)}
\newcommand{\qpt}[2]{\mathcal{E}_{{#1},{#2}}(x,w)}

\newcommand{\rxw}{\mathcal{R}(\protect{x},\protect{w})}
\newcommand{\rytw}{\mathcal{R}(yt,w)}
\newcommand{\rysw}{\mathcal{R}(ys,w)}
\newcommand{\rysiw}{\mathcal{R}(ys_i,w)}
\newcommand{\qabxw}{\mathcal{E}_{a,b}(x,w)}
\newcommand{\qptxw}{\mathcal{E}_t(x,w)}
\newcommand{\qptabxw}{\mathcal{E}_{\tab}(x,w)}
\newcommand{\qbcxw}{\mathcal{E}_{b,c}(x,w)}
\newcommand{\qcdxw}{\mathcal{E}_{c,d}(x,w)}
\newcommand{\qbcyw}{\mathcal{E}_{b,c}(y,w)}

\newcommand{\qabxtwt}{\mathcal{E}_{a,b}(\widetilde{x},\widetilde{w})}
\newcommand{\qptxtwt}{\mathcal{E}_t(\widetilde{x},\widetilde{w})}
\newcommand{\qptabxtwt}{\mathcal{E}_{\tab}(\widetilde{x},\widetilde{w})}
\newcommand{\qbcxtwt}{\mathcal{E}_{b,c}(\widetilde{x},\widetilde{w})}
\newcommand{\qcdxtwt}{\mathcal{E}_{c,d}(\widetilde{x},\widetilde{w})}

\newcommand{\rxiwi}{\mathcal{R}(x^{-1},w^{-1})}
\newcommand{\rysiwi}{\mathcal{R}(sy^{-1},w^{-1})}
\newcommand{\ryw}{\mathcal{R}(y,w)}
\newcommand{\ryiwi}{\mathcal{R}(y^{-1},w^{-1})}
\newcommand{\rxtwt}{\mathcal{R}(\widetilde{x},\widetilde{w})}
\newcommand{\rxtw}{\mathcal{R}(xt,w)}

\newcommand{\ykm}{y_{k,m}}
\newcommand{\xkm}{x_{k,m}}
\newcommand{\xktm}{y_{k,m}}
\newcommand{\xrs}{x_{r,s}}
\newcommand{\xrts}{y_{r,s}}
\newcommand{\wkm}{w_{k,m}}
\newcommand{\vkm}{v_{k,m}}
\newcommand{\wktm}{v_{k,m}}
\newcommand{\wrs}{w_{r,s}}
\newcommand{\wrts}{v_{r,s}}


%% file: jugmac.tex

\newcommand{\magnuss}{Magnus\ } 
\newcommand{\magnus}{Magnus} 

\newcommand{\emp}{\raisebox{-.2em}[-.2em][-.2em]{$\smallsmile$}}
\newcommand{\hemp}{\protect\raisebox{3pt}{$\smallsmile$}}
\newcommand{\cat}{\kern 1.385pt\bullet\kern 1.385pt}
\newcommand{\starr}[1]{\stackrel{#1}{\longrightarrow}}

\newcommand{\adv}[1]{\Theta_{#1}}
\newcommand{\shc}[1]{[#1\negmedspace\overset{\cat}{\leftharpoondown}]\,}
\newcommand{\she}[1]{[#1\negmedspace\leftharpoondown\kern -11.5pt \hemp]}
\newcommand{\invshc}[1]{[\,\stackrel{\cat}{\rightharpoondown}#1]}
\newcommand{\invshe}[1]{[\rightharpoondown\kern -11.5pt \hemp\,#1]}
\newcommand{\sminvshe}[1]{[\rightharpoondown\kern -8.5pt \hemp\,\,#1]}
\newcommand{\insc}[2]{I^{\cat}_{#1}(#2)}
\newcommand{\inse}[2]{I^{\emp}_{#1}(#2)}
\newcommand{\thr}{\tau}
\newcommand{\dro}{\delta}
\newcommand{\add}{\alpha}

\newcommand{\ltms}[1]{\operatorname{L-Times(#1)}}
\newcommand{\tms}[1]{\operatorname{Times(#1)}}
\newcommand{\snuo}{\nu_1}
\newcommand{\snur}{\nu_3}
\newcommand{\snui}{\nu_i}
\newcommand{\snut}{\nu_2}
\newcommand{\snuh}{\nu_h}
\newcommand{\snuhm}{\nu_{h-1}}
\newcommand{\snuj}{\nu_{j}}
\newcommand{\snujm}{\nu_{j-1}}
\newcommand{\snujp}{\nu_{j+1}}
\newcommand{\snujq}{\nu_{j+2}}

\newcommand{\nub}{\overline{\nu}}
\newcommand{\bnuo}{\nub_1}
\newcommand{\bnut}{\nub_2}
\newcommand{\bnur}{\nub_3}
\newcommand{\bnui}{\nub_i}
\newcommand{\bnuh}{\nub_h}
\newcommand{\bnuhm}{\nub_{h-1}}
\newcommand{\bnujm}{\nub_{j-1}}
\newcommand{\bnujp}{\nub_{j+1}}

\newcommand{\omo}{\omega_1}
\newcommand{\omr}{\omega_3}
\newcommand{\omi}{\omega_i}
\newcommand{\omt}{\omega_2}
\newcommand{\omh}{\omega_h}
\newcommand{\omhm}{\omega_{h-1}}
\newcommand{\omjm}{\omega_{j-1}}
\newcommand{\omjp}{\omega_{j+1}}

\newcommand{\wt}{\Delta}
\newcommand{\shdwt}{\Delta_\downarrow}
\newcommand{\shuwt}{\Delta_\uparrow}
\newcommand{\balph}{\boldsymbol{\alpha}}
\newcommand{\av}{\boldsymbol{a}}
\newcommand{\bv}{\boldsymbol{b}}
\newcommand{\ek}{\boldsymbol{e}^k}
\newcommand{\ep}{\boldsymbol{e}^p}
\newcommand{\freed}[1]{|#1|}
\newcommand{\nball}[1]{b(#1)}
\newcommand{\ind}[1]{{\av}_{#1}}

\newcommand{\sti}[2]{\genfrac{\{}{\}}{0pt}{}{#1}{#2}}
\newcommand{\bell}[1]{B_{#1}}

\newcommand{\stp}[1]{St_{#1}}
\newcommand{\stpb}[1]{\overline{St_{#1}}}
\newcommand{\stn}{St_h}
\newcommand{\stnb}{\widehat{St}_h}
\newcommand{\stnm}{St_{h-1}}
\newcommand{\stnmk}{St_{h-1,f}}
\newcommand{\stnmkm}{St_{h-1,f-1}}
\newcommand{\stnk}{St_{h,f}}
\newcommand{\stnkb}{\widehat{St}_{h,f}}
\newcommand{\stnkp}{St_{h,f+1}}
\newcommand{\stnpk}{St_{h+1,f}}

\newcommand{\edges}{\operatorname{Edges}}

\newcommand{\grastd}[2]{G_{#1,#2}}
\newcommand{\hgrastd}[2]{\widehat{G}_{#1,#2}}
\newcommand{\graddd}[1]{G_{#1}^{a-d}}
\newcommand{\grann}[1]{G_{#1}^{ann}}

\newcommand{\stw}{\boldsymbol{\omega}}
\newcommand{\stv}{\boldsymbol{\nu}}
\newcommand{\stu}{\boldsymbol{\mu}}
\newcommand{\stwr}{\omega}
\newcommand{\stvr}{\nu}
\newcommand{\stur}{\mu}
\newcommand{\hsnu}{\widehat{\boldsymbol{\nu}}}
\newcommand{\hsmu}{\widehat{\boldsymbol{\mu}}}
\newcommand{\hsw}{\widehat{\boldsymbol{\omega}}}
\newcommand{\hnu}{\widehat{\nu}}
\newcommand{\hmu}{\widehat{\mu}}

\newcommand{\btv}{\overline{\boldsymbol{\nu}}}
\newcommand{\btvr}{\overline{\nu}}
\newcommand{\btw}{\overline{\boldsymbol{\omega}}}

\newcommand{\mc}[2]{\operatorname{MC}(#1,#2)}

\setlength{\fboxsep}{0pt}
\newcolumntype{Z}{l@{\hspace{2pt}}}
\newcolumntype{Y}{l@{\hspace{1.8pt}}}

\newcommand{\markmat}{P}
\newcommand{\marktrans}[2]{p_{#1,#2}}
\newcommand{\probvec}{\boldsymbol{\Omega}}

\DeclareRobustCommand{\nicefrac}[3][\mathit]{\hspace{0.1em}%
  \raisebox{0.4ex}{$#1{\scriptstyle
#2}$}\hspace{-0.1em}/\hspace{-0.07em}%
  \mbox{$#1{\scriptstyle #3}$}}
\newcommand{\mf}[2]{\nicefrac[\mathrm]{#1}{#2}}

\newcommand{\fancyfrac}[2]{\raisebox{0.3ex}{#1}
\hspace*{-0.1em}/\hspace*{-0.1em}\raisebox{-0.3ex}{#2}}

\newcommand{\odoxdx}{{\,\fbox{\kern-3pt\begin{tabular}{ZZZZZZ} 
   $\emp$ & $\cdots$ & $\emp$ & $\cat$ & $\cdots$ & $\cat$ \end{tabular}}\,}}

\newcommand{\xooxox}{{\,\fbox{\kern-3pt\begin{tabular}{ZZZZZZ} 
   $\cat$ & $\emp$ & $\emp$ & $\cat$ & $\emp$ & $\cat$ \end{tabular}}\,}}
\newcommand{\oxxoxo}{{\,\fbox{\kern-3pt\begin{tabular}{ZZZZZZ}
   $\emp$ & $\cat$ & $\cat$ & $\emp$ & $\cat$ & $\emp$ \end{tabular}}\,}}
\newcommand{\xxoxoo}{{\,\fbox{\kern-3pt\begin{tabular}{ZZZZZZ}
   $\cat$ & $\cat$ & $\emp$ & $\cat$ & $\emp$ & $\emp$ \end{tabular}}\,}}
\newcommand{\xoxoox}{{\,\fbox{\kern-3pt\begin{tabular}{ZZZZZZ}
   $\cat$ & $\emp$ & $\cat$ & $\emp$ & $\emp$ & $\cat$ \end{tabular}}\,}}

\newcommand{\xxxxxxxoxxooxo}{{\,\fbox{\kern-3pt\begin{tabular}{ZZZZZZZZZZZZZZ}
  $\cat$ & $\cat$ & $\cat$ & $\cat$ & $\cat$ & 
  $\cat$ & $\cat$ & $\emp$ & $\cat$ & $\cat$ & 
  $\emp$ & $\emp$ & $\cat$ & $\emp$ \end{tabular}}\,}}

\newcommand{\xxxxxoxoxxooxo}{{\,\fbox{\kern-3pt\begin{tabular}{ZZZZZZZZZZZZZZ}
  $\cat$ & $\cat$ & $\cat$ & $\cat$ & $\cat$ & 
  $\emp$ & $\cat$ & $\emp$ & $\cat$ & $\cat$ & 
  $\emp$ & $\emp$ & $\cat$ & $\emp$ \end{tabular}}\,}}

\newcommand{\xxoxoxxxoo}{{\,\fbox{\kern-3pt\begin{tabular}{ZZZZZZZZZZ}
  $\cat$ & $\cat$ & $\emp$ & $\cat$ & $\emp$ & 
  $\cat$ & $\cat$ & $\cat$ & $\emp$ & $\emp$ \end{tabular}}\,}}

\newcommand{\aaoaoaaaoo}{{\,\fbox{\kern-3pt\begin{tabular}{ZZZZZZZZZZ}
  6 & 6 & $\emp$ & 6 & $\emp$ & 
  9 & 7 & 9 & $\emp$ & $\emp$ \end{tabular}}\,}}

\newcommand{\bbobobbboo}{{\,\fbox{\kern-3pt\begin{tabular}{ZZZZZZZZZZ}
  5 & 9 & $\emp$ & 9 & $\emp$ & 
  8 & 8 & 8 & $\emp$ & $\emp$ \end{tabular}}\,}}

\newcommand{\nuh}{{\,\fbox{\kern-3pt\begin{tabular}{ZZZZ}
   $\snuo$ & $\snut$ & $\cdots$ & $\snuh$ \end{tabular}}\,}}
\newcommand{\shenuh}{{\,\fbox{\kern-3pt\begin{tabular}{ZZZZ}
   $\snut$ & $\cdots$ & $\snuh$ & $\emp$ \end{tabular}}\,}}
\newcommand{\inshnuh}{{\,\fbox{\kern-3pt\begin{tabular}{ZZZZ}
   $\emp$ & $\snuo$ & $\cdots$ & $\snuhm$ \end{tabular}}\,}}
\newcommand{\inhnuh}{{\,\fbox{\kern-3pt\begin{tabular}{ZZZZ}
   $\emp$ & $\hnu_1$ & $\cdots$ & $\hnu_{h-1}$ \end{tabular}}\,}}
\newcommand{\insjn}{{\,\fbox{\kern-3pt\begin{tabular}{ZZZZZZZZ}
   $\snuo$ & $\snut$ & $\cdots$ & $\snujm$ & $\emp$ & 
   $\snujp$ & $\cdots$ & $\snuh$ \end{tabular}}\,}}

\newcommand{\insjq}{{\,\fbox{\kern-3pt\begin{tabular}{YYYYYYYY}
   $\snut$ & $\cdots$ & $\snuj$ & $\cat$ & 
   $\snujq$ & $\cdots$ & $\snuh$ & $\emp$ \end{tabular}}\,}}

\newcommand{\xxxoo}{{\,\fbox{\kern-3pt\begin{tabular}{ZZZZZ}
   $\cat$ & $\cat$ & $\cat$ & $\emp$ & $\emp$ \end{tabular}}\,}}
\newcommand{\xxoxo}{{\,\fbox{\kern-3pt\begin{tabular}{ZZZZZ}
   $\cat$ & $\cat$ & $\emp$ & $\cat$ & $\emp$ \end{tabular}}\,}}
\newcommand{\xxoox}{{\,\fbox{\kern-3pt\begin{tabular}{ZZZZZ}
   $\cat$ & $\cat$ & $\emp$ & $\emp$ & $\cat$ \end{tabular}}\,}}
\newcommand{\xoxxo}{{\,\fbox{\kern-3pt\begin{tabular}{ZZZZZ}
   $\cat$ & $\emp$ & $\cat$ & $\cat$ & $\emp$ \end{tabular}}\,}}
\newcommand{\xoxox}{{\,\fbox{\kern-3pt\begin{tabular}{ZZZZZ}
   $\cat$ & $\emp$ & $\cat$ & $\emp$ & $\cat$ \end{tabular}}\,}}
\newcommand{\xooxx}{{\,\fbox{\kern-3pt\begin{tabular}{ZZZZZ}
   $\cat$ & $\emp$ & $\emp$ & $\cat$ & $\cat$ \end{tabular}}\,}}
\newcommand{\oxxxo}{{\,\fbox{\kern-3pt\begin{tabular}{ZZZZZ}
   $\emp$ & $\cat$ & $\cat$ & $\cat$ & $\emp$ \end{tabular}}\,}}
\newcommand{\oxxox}{{\,\fbox{\kern-3pt\begin{tabular}{ZZZZZ}
   $\emp$ & $\cat$ & $\cat$ & $\emp$ & $\cat$ \end{tabular}}\,}}
\newcommand{\oxoxx}{{\,\fbox{\kern-3pt\begin{tabular}{ZZZZZ}
   $\emp$ & $\cat$ & $\emp$ & $\cat$ & $\cat$ \end{tabular}}\,}}
\newcommand{\ooxxx}{{\,\fbox{\kern-3pt\begin{tabular}{ZZZZZ}
   $\emp$ & $\emp$ & $\cat$ & $\cat$ & $\cat$ \end{tabular}}\,}}

\newcommand{\oxxxoooxoxx}{{\,\fbox{\kern-3pt \begin{tabular}{ZZZZZZZZZZZ}
   $\emp$ & $\cat$ & $\cat$ & $\cat$ & $\emp$ & $\emp$ & 
   $\emp$ & $\cat$ & $\emp$ & $\cat$ & $\cat$ \end{tabular}}\,}}

\newcommand{\xxxo}{{\,\fbox{\kern-1pt\begin{tabular}{ZZZZ}
        $\cat$ & $\cat$ & $\cat$ & $\emp$ \end{tabular}}}}
                                                                              
\newcommand{\xxox}{{\,\fbox{\kern-1pt\begin{tabular}{ZZZZ}
        $\cat$ & $\cat$ & $\emp$ & $\cat$ \end{tabular}}}}
                                                                              
\newcommand{\xoxx}{{\,\fbox{\kern-1pt\begin{tabular}{ZZZZ}
        $\cat$ & $\emp$ & $\cat$ & $\cat$ \end{tabular}}}}
                                                                              
\newcommand{\oxxx}{{\,\fbox{\kern-1pt\begin{tabular}{ZZZZ}
        $\emp$ & $\cat$ & $\cat$ & $\cat$ \end{tabular}}}}


\psfrag{6-46--7}{  \fbox{\kern-1pt\begin{tabular}{ZZZZZZZ}$6$ & $\emp$ & $4$ &
      $6$ & $\emp$ & $\emp$ & $7$\end{tabular}}}

\psfrag{-------}{  \fbox{\kern-1pt\begin{tabular}{ZZZZZZZ}$\emp$ & $\emp$ & $\emp$ &
      $\emp$ & $\emp$ & $\emp$ & $\emp$\end{tabular}}}


\psfrag{brackxx-}{  \raisebox{-1.6pt}{\huge [}\,\fbox{\kern-1pt\begin{tabular}{ZZZ}$\cat$ & $\cat$ & $\emp$ \end{tabular}} \raisebox{-1.6pt}{\huge ]}}
                                                                              
\psfrag{brackx-x}{ \raisebox{-1.6pt}{\huge [}\,\fbox{\kern-1pt\begin{tabular}{ZZZ}$\cat$ & $\emp$ & $\cat$ \end{tabular}} \raisebox{-1.6pt}{\huge ]}}
                                                                              
\psfrag{brack-xx}{ \raisebox{-1.6pt}{\huge [}\,\fbox{\kern-1pt\begin{tabular}{ZZZ}$\emp$ & $\cat$ & $\cat$ \end{tabular}} \raisebox{-1.6pt}{\huge ]}}
                                                                              
\psfrag{ac-}{ \fbox{\kern-1pt\begin{tabular}{ZZZ}$1$ & $3$ & $\emp$ \end{tabular}}}

\psfrag{cb-}{ \fbox{\kern-1pt\begin{tabular}{ZZZ}$3$ & $2$ & $\emp$ \end{tabular}}}

\psfrag{bb-}{ \fbox{\kern-1pt\begin{tabular}{ZZZ}$2$ & $2$ & $\emp$ \end{tabular}}}

\psfrag{cc-}{ \fbox{\kern-1pt\begin{tabular}{ZZZ}$3$ & $3$ & $\emp$ \end{tabular}}}

\psfrag{b-c}{ \fbox{\kern-1pt\begin{tabular}{ZZZ}$2$ & $\emp$ & $3$ \end{tabular}}}

\psfrag{c-c}{ \fbox{\kern-1pt\begin{tabular}{ZZZ}$3$ & $\emp$ & $3$ \end{tabular}}}

\psfrag{-cc}{ \fbox{\kern-1pt\begin{tabular}{ZZZ}$\emp$ & $3$ & $3$ \end{tabular}}}


\psfrag{xxx-}{ \fbox{\kern-1pt\begin{tabular}{ZZZZ}$\cat$ & $\cat$ & $\cat$ & $\emp$ \end{tabular}}}
                                                                              
\psfrag{xx-x}{ \fbox{\kern-1pt\begin{tabular}{ZZZZ}$\cat$ & $\cat$ & $\emp$ & $\cat$ \end{tabular}}}
                                                                              
\psfrag{x-xx}{ \fbox{\kern-1pt\begin{tabular}{ZZZZ}$\cat$ & $\emp$ & $\cat$ & $\cat$ \end{tabular}}}
                                                                              
\psfrag{-xxx}{ \fbox{\kern-1pt\begin{tabular}{ZZZZ}$\emp$ & $\cat$ & $\cat$ & $\cat$ \end{tabular}}}

\psfrag{xxx--}{ \fbox{\begin{tabular}{ZZZZZ} $\cat$ & $\cat$ & $\cat$ & $\emp$ & $\emp$ \end{tabular}}}
\psfrag{xx-x-}{ \fbox{\begin{tabular}{ZZZZZ} $\cat$ & $\cat$ & $\emp$ & $\cat$ & $\emp$ \end{tabular}}}
\psfrag{xx--x}{ \fbox{\begin{tabular}{ZZZZZ} $\cat$ & $\cat$ & $\emp$ & $\emp$ & $\cat$ \end{tabular}}}
\psfrag{x-xx-}{ \fbox{\begin{tabular}{ZZZZZ} $\cat$ & $\emp$ & $\cat$ & $\cat$ & $\emp$ \end{tabular}}}
\psfrag{x-x-x}{ \fbox{\begin{tabular}{ZZZZZ} $\cat$ & $\emp$ & $\cat$ & $\emp$ & $\cat$ \end{tabular}}}
\psfrag{x--xx}{ \fbox{\begin{tabular}{ZZZZZ} $\cat$ & $\emp$ & $\emp$ & $\cat$ & $\cat$ \end{tabular}}}
\psfrag{-xxx-}{ \fbox{\begin{tabular}{ZZZZZ} $\emp$ & $\cat$ & $\cat$ & $\cat$ & $\emp$ \end{tabular}}}
\psfrag{-xx-x}{ \fbox{\begin{tabular}{ZZZZZ} $\emp$ & $\cat$ & $\cat$ & $\emp$ & $\cat$ \end{tabular}}}
\psfrag{-x-xx}{ \fbox{\begin{tabular}{ZZZZZ} $\emp$ & $\cat$ & $\emp$ & $\cat$ & $\cat$ \end{tabular}}}
\psfrag{--xxx}{ \fbox{\begin{tabular}{ZZZZZ} $\emp$ & $\emp$ & $\cat$ & $\cat$ & $\cat$ \end{tabular}}}

\psfrag{---}{ \fbox{\begin{tabular}{ZZZ} $\emp$ & $\emp$ & $\emp$ \end{tabular}}}
\psfrag{x--}{ \fbox{\begin{tabular}{ZZZ} $\cat$ & $\emp$ & $\emp$ \end{tabular}}}
\psfrag{-x-}{ \fbox{\begin{tabular}{ZZZ} $\emp$ & $\cat$ & $\emp$ \end{tabular}}}
\psfrag{--x}{ \fbox{\begin{tabular}{ZZZ} $\emp$ & $\emp$ & $\cat$ \end{tabular}}}
\psfrag{xx-}{ \fbox{\begin{tabular}{ZZZ} $\cat$ & $\cat$ & $\emp$ \end{tabular}}}
\psfrag{x-x}{ \fbox{\begin{tabular}{ZZZ} $\cat$ & $\emp$ & $\cat$ \end{tabular}}}
\psfrag{-xx}{ \fbox{\begin{tabular}{ZZZ} $\emp$ & $\cat$ & $\cat$ \end{tabular}}}
\psfrag{xxx}{ \fbox{\begin{tabular}{ZZZ} $\cat$ & $\cat$ & $\cat$ \end{tabular}}}

